
\documentstyle{amsppt}
\magnification=1200
\hsize=5.85truein
\vsize=8.15truein
\voffset 0.25truein
\TagsOnRight
\baselineskip=14pt
\parindent 1.5em
\nopagenumbers

\document

\rightline {ALBERTA THY 32-95}

\bigskip

\bigskip

\bigskip

\centerline {\bf GRADED CONTRACTIONS OF AFFINE KAC-MOODY ALGEBRAS\/}

\bigskip
\bigskip
\centerline {Marc de Montigny}

\centerline {Facult\'e Saint-Jean and Theoretical Physics Institute}

\centerline {University of Alberta, 8406 91 Street}

\centerline {Edmonton, Alberta, T6C 4G9}

\centerline {CANADA}

\bigskip

\bigskip

\bigskip

\centerline {\bf December 1995\/}

\bigskip

\bigskip

\bigskip

ABSTRACT.
The method of
 graded contractions,
 based on the
 preservation of the automorphisms of finite order,
 is applied to the affine Kac-Moody algebras and their
 representations, to yield a new class
 of infinite dimensional Lie algebras and representations.
 After the introduction of the horizontal and
 vertical gradings, and the algorithm to find the horizontal toroidal
 gradings, I discuss some general properties of the
 graded contractions, and compare them with the
 In\"on\"u-Wigner contractions.
 The example of $\hat A_2$ is discussed in detail.

\bigskip

\leftline {{\it PACS numbers:} 02.20.+b, 11.30.--j}

\bigskip

\leftline {{\it Keywords:} Lie algebras, representation theory, contraction.}

\bigskip

\newpage

\leftline {\bf 1. INTRODUCTION.}
\medskip

In this paper, I
 describe the {\it graded
 contractions} of general affine Kac-Moody algebras and
 their representations, and illustrate the method
 with $\hat A_2$. Generally speaking,
 contractions of Lie algebras are deformations, or singular
 transformations, of the constants of structure. They were
 introduced in physics by E. In\"on\"u and E. Wigner
 \cite{1} in order to provide
 a formal relationship between the kinematical
 groups of Einstein's special relativity and Galilean
 relativity.
 In general, contractions are interesting because they relate, in
 a meaningful way, different Lie algebras such that various properties
 of the contracted (or limit) algebra can be obtained from the initial
 algebra. This is
 particularly promising in the study of non-semisimple Lie
 algebras (which often can be seen as
 the outcome of a contraction procedure),
 because their representation theory and their general structure
 are not as elegant and uniform as the semisimple Lie
 algebras.

Although the Lie algebras most familiar to physicists are
 {\it simple}, such as $su(2)$, $su(3)$ or $E_6$, many algebras of physical
 interest are likely to be non-semisimple. The situation is similar with the
 infinite dimensional Lie algebras. Since the Kac-Moody and
 Virasoro algebras represent a rather restricted class of algebras,
 their contractions lead to a totally new class of infinite
 dimensional Lie algebras, which might as well be relevant in physics.
 An example of an infinite dimensional algebra which can be seen as
 a contraction of a Kac-Moody algebra is the
 oscillator (or Heisenberg) algebra, with commutation relations:
$$
 [\hslash , a_n]=0,\qquad [a_m,a_n]=m\delta_{m+n,0}\hslash .
$$
 Since the early work of In\"on\"u and Wigner,
 the method has been generalized in many directions
 (some of which are given in  \cite{2--7}) and
 have been applied to various problems in physics (see, for instance,
  \cite{8--15}).
 To my knowledge, the first systematical
 treatment of In\"on\"u-Wigner contractions of Kac-Moody algebras appeared in
  \cite{16}.

Similar ``limits'' of infinite dimensional Lie algebras already
 exist (implicitly) in the literature. For instance, in their
 investigation of the principal chiral field model in low dimensions,
 Faddeev and Reshetikhin \cite{13}
 have considered the current algebra
$$
\aligned
&[S^a(x), S^b(y)]=i\epsilon^{abc}S^c(x)\delta (x-y)-
\dfrac {ik}{2\pi}\delta^{ab}\delta^\prime(x-y),\\
&[T^a(x),T^b(y)]=i\epsilon^{abc}T^c(x)\delta(x-y)+
\dfrac {ik}{2\pi}\delta^{ab}\delta^\prime (x-y),\\
&[S^a(x),T^b(y)]=0,\endaligned$$
which, under the change of basis $A(x)\equiv \gamma(S+T)(x)$ and
 $B(x)\equiv\dfrac\pi k(S-T)(x)$, becomes, as $k$ approaches infinity,
$$
\aligned
&[A^a(x), A^b(y)]=i\epsilon^{abc}\gamma A^c(x)\delta (x-y),\\
&[A^a(x), B^b(y)]=i\epsilon^{abc}\gamma B^c(x)\delta (x-y)
 -i\gamma\delta^{ab}\delta^\prime (x-y),\\
&[B^a(x), B^b(y)]=0.\endaligned$$

Recently, non-semisimple affine Kac-Moody algebras have been used
 in string theory, in the context of Wess-Zumino-Witten
 (WZW) models \cite{17}, where they occur in the
 expression of the current algebras of the models. String
 backgrounds based on non-semisimple WZW models have been constructed in
 \cite{18} and a general class of exact conformal field theories,
 with integral Virasoro central charges, have been constructed in
 \cite{19}. Other constructions appear in  \cite{20}.
 Although these models have been the main motivation for the present work,
 it has interest of its own and is not restricted to these
 applications.

In this paper, I
 apply a method of contraction \cite{22,23}
 based on the preservation
 of a {\it grading}: a decompositon of the Lie algebra into eigenspaces of
 an automorphism of finite order. I describe the algorithm to find
 the grading preserving contractions, or {\it graded contractions},
 of an affine Kac-Moody algebra. Starting with an affine algebra
 we obtain different ({\it i.e.} non-isomorphic) infinite
 dimensional Lie algebras. The interest of this particular method is
 when we require, for some physical reasons, one or many
 automorphisms ({\it e.g.} parity or time reversal) to
 be admitted by the limit algebras. The systematical study of Lie gradings
 has been initiated in \cite {21}, as a powerful tool in Lie
 theory. When it comes to non-semisimple algebras, the graded
 contractions could be most useful in studying the gradings.
 Indeed, in some
 cases, a contraction is the only way to build representations of ``exotic''
 infinite algebras, whose representation theory is yet to be
 understood. To summarize, whenever one has to use
 some properties of a non-semisimple Lie algebras,  many of
 these properties can be obtained from a contraction, and if
 an automorphism of finite order is preserved through that contraction,
 then the formalism of graded contractions is possibly more appropriate.
 An advantage of this method is that it applies simultaneously to all
 the algebras and representations which admit a common grading
 structure.

 The method described in \cite {22, 23} is
 implicitly applicable to infinite dimensional Lie algebras, but
 it is studied systematically
 (for all affine Kac-Moody algebras) for the first time here.
 The particular case $\hat A_1$ has been considered in \cite{25}.
 I
 also introduce the concept of
 {\it vertical} and {\it horizontal} gradings, which do not exist
 with the finite dimensional algebras. With the purpose of applying
 these results to high energy physics, I consider also the contractions
 of (integrable irreducible highest weight) representations, and
 discuss various properties of the contracted algebras.
 Although the present formalism
 can be applied to superalgebras, I do not consider them
 here.

I close this section by reviewing briefly the method of graded
 contractions, introduced in \cite {22, 23} and reviewed in
 \cite{24}.

\medskip

\leftline {\it Definition of graded contractions.}
\smallskip

A {\it grading} of a (finite or infinite dimensional)
 Lie algebra ${\frak g}$
 is a vector decomposition:
$$
\frak g=\bigoplus_{\mu\in\Gamma}\frak g_\mu,
\quad\text {such that}\quad \left [{\frak g}_\mu,{\frak g}_\nu\right ]\subseteq
{\frak g}_{\mu +\nu},\tag1.1$$
where $\mu$ and $\nu$ belong to a finite abelian grading group
 $\Gamma =\Bbb Z_{m1}\otimes\cdots\otimes\Bbb Z_{mk}$,
 and the notation
 $\left [{\frak g}_\mu,{\frak g}_\nu\right ]$ means that if
 $x\in{\frak g}_\mu$ and  $y\in{\frak g}_\nu$, then either
 $[x,y]\subseteq{\frak g}_{\mu +\nu}$, or it is zero
 (${\frak g}_\mu, {\frak g}_nu$
 represent any element of the respective grading subspaces).
 The product in the grading group $\Gamma$ is denoted
 by $+$.

Along with the decomposition (1.1), a grading of a $\frak
 g$-module $V$ is a splitting:
$$
V=\bigoplus_{\mu\in\Gamma} V_\mu,
\quad \text {with}\
\frak g_\mu \cdot V_\nu\subseteq V_{\mu+\nu},\tag1.2$$
where the action $\sigma_V(\frak g)V$ is denoted
 $\frak g\cdot V$. The expression on the right-hand side has a meaning
 similar to (1.1).
 As mentioned previously, a grading is associated to
 an element $\phi\in Aut(\frak g)$ of
 finite order ({\it i.e.} $\phi^M=id_{Aut(\frak g)}$, $M$ finite)
 that acts on $\frak g$ and $V$ as
$
\phi g=\left (\text{exp}\frac {2\pi i}Mk\right )g,$
if $g\in{\frak g}_k$, and
$
\phi v=\left (\text{exp}\frac {2\pi i}Ml\right )v,$
if $v\in V_l$. In general, this automorphism comes from a physical
 restriction. For instance, in \cite {20} one may notice that, in the
 construction of a WZW model, if it is possible to split the initial
 group into the
 ``coset''part and the ``subgroup''part by the matrix that rotate
 the generators, then it follows that
 if this matrix is an automorphism, then the
 contracted algebra admits this automorphism as well.

The contractions which preserve the $\Gamma$ gradings
 (1.1), (1.2) are called {\it graded contractions}.  The graded
 contractions of the Lie algebra $\frak g$ are defined by introducing
 parameters $\varepsilon_{\mu,\nu}$, such that the contracted
 algebra ${\frak g}^\varepsilon$ has the same vector space
 as the original algebra $\frak g$, but
 modified commutation relations:
$$
\left [{\frak g}_\mu,{\frak g}_\nu\right ]_\varepsilon\equiv
 \varepsilon_{\mu,\nu}\left [{\frak g}_\mu,{\frak g}_\nu\right ]\subseteq
 \varepsilon_{\mu,\nu}{\frak g}_{\mu +\nu}.\tag 1.3$$
Similarly, the graded contractions of representations
 are defined through the introduction of parameters $\psi_{\mu,\nu}$,
 which deforms the action of $\frak g$ on $V$ such that they preserve
 the grading (1.2):
$$
{\frak g}_\mu{\cdot^\psi} V_\nu\equiv\psi_{\mu,\nu}{\frak g}_\mu\cdot
V_\nu \subseteq\psi_{\mu,\nu}V_{\mu+\nu}.\tag1.4$$
In \cite{22, 23}, it is shown that the
 contraction parameters $\varepsilon$
 and $\psi$ must satisfy the equations
$$
\varepsilon_{\mu,\nu}\varepsilon_{\mu+\nu,\sigma}=
\varepsilon_{\nu,\sigma}\varepsilon_{\sigma+\nu,\mu}=
\varepsilon_{\mu,\sigma}\varepsilon_{\mu+\sigma,\nu},\tag1.5$$
and
$$
\varepsilon_{\mu,\nu}\psi_{\mu+\nu,\sigma}=
\psi_{\nu,\sigma}\psi_{\mu,\nu+\sigma}=
\psi_{\mu,\sigma}\psi_{\nu,\mu+\sigma}.\tag1.6$$
The solutions of these two sets of equations,
 substituted back into (1.3) and (1.4), provide the contractions of
 the algebra ${\frak g}$ and its representations.
 To each set of parameters $\varepsilon$
 (which defines a contracted algebra), the corresponding
 solutions of (1.6) for the $\psi$'s yield contractions of the
 representation $V$. More details and remarks
 are given in \cite {22--24}.

Finally, let me mention that if $V$ and $W$ are two
 compatibly $\Gamma$ graded
 $\frak g$-modules, then their tensor product space $V\otimes W$ is
 graded according to
$$
V\otimes W=\bigoplus_{\sigma\in\Gamma} (V\otimes W)_\sigma,
 \quad\text {where}\ \
(V\otimes W)_\sigma\equiv\bigoplus_{\mu+\nu=\sigma}V_\mu\otimes W_\nu.
\tag1.7$$
The tensor product is deformed by introducing further
 contraction parameters $\tau$:
$$
(V\otimes W)^\tau_\sigma\equiv\bigoplus_{\mu+\nu}
\tau_{\mu,\nu}V_\mu\otimes W_\nu,\tag1.8$$
and the (symmetric) $\tau$ parameter must satisfy \cite{23,24}
$$
\psi_{\sigma,\mu+\nu}\tau_{\mu,\nu}=
\psi_{\sigma,\mu} \tau_{\sigma+\mu,\nu}=
\psi_{\sigma,\nu} \tau_{\sigma+\nu,\nu}.\tag1.9$$

In the next section, I present the concepts of vertical and
 horizontal gradings of Kac-Moody algebras, and illustrate them with
 $\hat A_2$. In section 3 I discuss the graded contractions
 and some properties.
 The purpose of this paper is not to
 provide huge (and not particularly useful)
 tables of contracted algebras.
 Instead, one can rely on the program presented
 in \cite{24}, given a specific grading or algebra.
 I rather describe gradings of Kac-Moody algebras and
 some general features of the contractions, emphasizing the main
 differences with the traditional methods.
 In fact, the less trivial part of the method is always to
 find the gradings (1.1-2), and then one just have to substitute
 the solutions $\varepsilon$ or $\psi$ into (1.3-4).
 There exists no uniform prescription to find all the gradings of a
 general Lie algebra.
 A comprehensive list of gradings exists only for
 the simple Lie algebras of rank two and some of rank three.
 However, I sketch in section 2.1 a method
 which provides an important class of gradings of semi-simple Lie
 algebras: the {\it toroidal gradings}.

\bigskip

\leftline {\bf 2. GRADINGS OF AFFINE KAC-MOODY ALGEBRAS.}
\medskip

Consider a simple complex Lie algebra $\frak g$, and the corresponding
 affine untwisted Kac-Moody
 algebra $\hat\frak g=\left (\frak g\otimes\Bbb C
 \left [t,t^{-1}\right ]\right )\oplus \Bbb C k$, where
 $\Bbb C \left [t,t^{-1}\right ]$ is the associative algebra of
 the Laurent polynomials in $t$, and
 $k$ is a central
 extension. The first term, $\tilde\frak g =\frak g\otimes\Bbb C
 \left [t,t^{-1}\right ]$, in the direct sum is called the
 {\it loop algebra} of ${\frak g}$.
 Given $a, b\in\frak g$,
 the commutation relations in $\hat\frak g$ read
$$
[a\otimes t^m, b\otimes t^n]=[a,b]\otimes t^{m+n} +
 m k B(a,b)\delta_{m+n,0},\tag2.1$$
where $m, n\in\Bbb Z$, $[a,b]$ is the commutator in ${\frak g}$,
 and $B(a,b)$ is the Killing form of $\frak g$. A few
 comments about the gradings of
 {\it twisted} algebras appear at the end of section
 2.2. General properties of infinite dimensional
 Lie algebras can be found in \cite {26--29} and references
 therein.

Here, I distinguish two classes of gradings:
 {\it horizontal} and {\it vertical}.
 The first is a grading of the finite algebra $\frak g$
 that is preserved through the affinization process. The
 vertical gradings are {\it not} inherent in $\frak g$, but are rather
 given by the gradings of $\Bbb C\left [t,t^{-1}\right ]$. Obviously,
 these two types of gradings can be combined to provide gradings
 which are neither vertical nor horizontal.
 One can think of these two
 types of gradings
 as being the building blocks of the gradings of $\hat\frak g$.
 Vertical gradings have no analogue in the finite
 dimensional algebras.

\medskip

\leftline {\it 2.1. Horizontal gradings.}
\smallskip

The gradings of a {\it finite} Lie algebra $\frak g$ are associated to
 automorphisms of finite order (or conjugacy classes of elements of
 finite order (EFO)) of the corresponding compact Lie group $K$
 \cite{30, 31}. (An elementary introduction to the EFO theory,
 sufficient for our purposes,
 is given in \cite{33}).
 Kac's theory of EFO provides a prescription
 to identify the conjugacy classes of EFO, and hence the gradings.
 Such as described below, the action of an EFO leads to a {\it toroidal}
 $\Bbb Z_N$ grading (for an EFO
 of order $N$) and one can use it to
 grade simultaneously a Lie algebra and its irreducible representations.
 (It is called ``toroidal'' because it is a coarsening
 of the {\it toroidal}
 --or {\it Cartan}-- decomposition).
 This method provides, in a straightforward way, the unique diagonal
 representative of a conjugacy class of EFO in any irreducible
 representation of $\frak g$. All one needs to know is the weight
 system of the representation.

To grade an algebra, one must
 consider the {\it adjoint} representation,
 for which the weight system is the
 root system of $\frak g$. If $\frak g$ has rank $r$, then the EFO
 is represented by
 an array of non-negative integers
 $[s_0,\dots ,s_r]$, with $1$ as a common
 divisor. To each root $\alpha$ --and therefore each basis element of
 $\frak g$-- is associated the eigenvalue
$$
\text {exp}\frac{2\pi i}M<\alpha,{\bold s}>,\tag2.2$$
where $<\alpha,{\bold s}>=\sum_{j=1}^ra_js_j$, if $\alpha$ is given by
 $\alpha=\sum_{j=1}^r a_j\alpha_j$ ($\alpha_j$ are the
 simple roots of $\frak g$). Also, $M=s_0+\sum_{j=1}^r c_js_j$, the
 $c_j$ being the components (called {\it marks}) of the highest root
 $\Psi=\sum_{j=1}^r c_j\alpha_j$ of $\frak g$. The vector of marks
 is annihilated by the affine Cartan matrix $A$ ({\it i.e.}
 $\sum_{k=0}^r c_jA_{jk}=0,$ for all $k$). In a Cartan-Weyl basis,
 the elements $e_\alpha$ belong to the eigenvalue given by (2.2), for
 any positive or negative root $\alpha$, and all the elements of the
 Cartan subalgebra obviously belong to the eigenvalue $1$.
 The {\it order} of the EFO is $N=MC$, where $C$ is given in Table 6
 of \cite{31} for all the simple Lie algebras. The grading group
 is then $\Gamma=\Bbb Z_N$.

This can be generalized to any weight system.
 Let $V(\Lambda )$ be an irreducible $\frak g$-module with highest
 weight $\Lambda$, that can be {\it Cartan}-decomposed as:
$$
V(\Lambda )=\bigoplus_{\lambda\in\Omega(\Lambda)}V^\Lambda_\lambda ,
\qquad V^\Lambda_\lambda =\{v\in V(\Lambda )|\ hv=\lambda (h)v\},\tag2.3$$
where $h$ is an element of $\frak h$, the Cartan subalgebra of
 $\frak g$, and $\Omega(\Lambda)$ is the weight system of the module.
 To each weight $\lambda\in\Omega(\Lambda)$ a grading decomposition
 (1.1) is obtained by determining eigenvalues similar to (2.2):
$$
v_\lambda\longrightarrow v_\lambda\ \text{exp}
 \frac{2\pi i}M<\lambda,{\bold s}>,\tag2.4$$
where $<\lambda,{\bold s}>=\sum_{j=1}^n b_js_j$, if
 $\lambda=\sum_{j=1}^nb_j\alpha_j$. The value of $M$ is the same as
 in (2.2). Obviously, there are other --{\it non-toroidal}-- gradings
 of ${\frak g}$ which can serve as horizontal gradings, but I do not
 consider them hereafter because they are often related
 to toroidal gradings.
 For example, a grading of $\frak g$ can be provided by an EFO of a
 larger group, which contains $\frak g$.
 A classification of such gradings does not exist yet.

A general $\Bbb Z_{m1}\otimes\cdots\otimes\Bbb Z_{mk}$ is obtained by
 ``mixing'' gradings $\Bbb Z_{m1}, \Bbb Z_{m2}, \dots$ found by
 using the EFO. If each $\Bbb Z_{mj}$ provides
 a decomposition of $\frak g$:
$$
{\frak g}=\bigoplus_{\mu_j\in\Bbb Z_{m_j}}{\frak g}_{\mu_j},\tag2.5$$
for $j=1,\dots k$, then a
 $\Gamma=\Bbb Z_{m1}\otimes\cdots\otimes\Bbb Z_{mk}$ grading is obtained
as follows:
$$
{\frak g}=\bigoplus_{\mu\in\Gamma}{\frak
g}_{\mu=(\mu_1\cdots\mu_k)},\tag2.6$$
where ${\frak g}_{(\mu_1\cdots\mu_k)}={\frak g}_{\mu_1}\cap
\cdots\cap{\frak g}_{\mu_k}$. An example
 is given at the end of the
 section 2.3.

Once we have a $\Gamma$ grading (1.1) of a {\it finite}
 Lie algebra ${\frak g}$ (obtained from the EFO, or
 otherwise), then the
 corresponding affine Lie algebra ${\hat\frak g}$ admits the
 horizontal grading:
$$
{\hat\frak g}=\bigoplus_{\mu\in\Gamma}{\hat\frak g}_\mu,
\tag2.7$$
where
$$
\aligned
&{\hat\frak g}_0=\left ({\frak g}_0\otimes
 \Bbb C\left [t,t^{-1}\right ]\right )\oplus \Bbb Ck,\\
&{\hat\frak g}_{\mu\neq 0}={\frak g}_\mu\otimes
 \Bbb C\left [t,t^{-1}\right ].\endaligned\tag2.8$$
The identity element of $\Gamma$ is denoted by $0$.
 Note that, whenever $\hat\frak g$ is
 contracted or not, each ${\hat\frak
g}_{\mu\neq 0}$ carries a representation space for the subalgebra
 ${\hat\frak g}_0$.

In order to find gradings of an irreducible integrable highest weight
 representation of ${\hat\frak g}$, it is useful
 to express the gradings of ${\hat\frak g}$ in terms of the
 root vectors of ${\hat\frak g}$. If $\alpha\in\Delta_{\frak g}$
 is a root in ${\frak g}$
 (with components $\alpha_1,\dots,\alpha_r$),
 then the element
 $e_\alpha\otimes t^m$ can be denoted by
 $E_{\alpha+m\delta}$.
 The root $\alpha_0$ is given by
 $\delta=\alpha_0+\Psi\ (\Psi$ is the highest root of ${\frak g}$).
 Therefore the root
 vector can be expressed solely in terms of the affine simple roots
 ($\alpha_0,\alpha_1,\dots ,\alpha_r$), or
 in terms of $\delta$ and the finite simple roots:
 ($\delta,\alpha_1,\dots,\alpha_r$).
 To any element $h_k$ in the
 Cartan subalgebra of ${\frak g}$ is associated the root vector
 $E^k_{m\delta}$. In the case of a $\Bbb Z_N$ grading provided by
 (2.2), $E_{\alpha+m\delta}$ belongs to the subspace
 ${\hat\frak g}_\mu$, for any $m$, if $e_\alpha$ belongs to the grading
 subspace ${\frak g}_\mu$. Obviously, all the vectors $E_{m\delta}$
 belong to ${\hat\frak g}_0$, as do the elements of the Cartan
 subalgebra of ${\hat\frak g}$.

Finally, I repeat that the gradings of a finite algebra are not
 always manifestly related to an EFO, and that there are other types
 of gradings
 ({\it e.g.} the generalized Pauli matrices used
 in \cite {21}). Such gradings are the result of an {\it outer}
 automorphism of ${\frak g}$, whereas the EFO correspond
 to {\it inner} automorphisms.
 In any event, once a grading (2.6) of a finite
 algebra is known,
 the equations (2.7)-(2.8) provide the corresponding horizontal grading
 of the affine algebra.

\medskip

\leftline {\it 2.2. Vertical gradings.}
\smallskip

Similarly, the vertical $\Bbb Z_N$ gradings are
 given by the action of a root of the unity
 $\exp\left (\dfrac {2\pi i}N\right )$
 on the associative algebra $\Bbb C\left [t,t^{-1}
 \right ]$:
$$
\phi: t\rightarrow \exp\left (\dfrac {2\pi i}N\right )\ t, \tag2.9$$
 such that the element $t^m$ belongs to the eigenvalue
 $\exp\left (\dfrac {2\pi i}Nm\right )$, and we write
 $t^m\in\Bbb C\left [t,t^{-1}\right ]_{m\ mod\ N}$. Therefore
 the grading can be written
$$
\Bbb C\left [t,t^{-1}\right ]=
\bigoplus_{j=0}^{N-1}\Bbb C\left [t,t^{-1}\right ]_j,\tag2.10$$
 where
$$
\Bbb C\left [t,t^{-1}\right ]_j=\oplus \Bbb C\ t^{j+kN},
 \qquad k\in\Bbb Z.\tag2.11$$

Accordingly, the grading of the Kac-Moody algebra
 ${\hat\frak g}$ is
$${\hat\frak g}=\bigoplus_{j\in\Bbb Z_N}{\hat\frak g}_j,\tag2.12$$
where
$$
\aligned
&{\hat\frak g}_0=({\frak g}\otimes\Bbb C\left [t,t^{-1}\right ]_0)
 \oplus\Bbb C k,\\
&{\hat\frak g}_j={\frak g}\otimes\Bbb C\left [t,t^{-1}\right ]_j.\endaligned
\tag2.8'$$

For example, a $\Bbb Z_3$ grading gives the decomposition:
$$
\aligned
&{\hat\frak g}_0=\Bbb C k+\cdots + {\frak g}\otimes t^{-3} +
{\frak g}\otimes t^0+{\frak g}\otimes t^3
+\cdots\\
&{\hat\frak g}_1= \cdots + {\frak g}\otimes t^{-2} +
{\frak g}\otimes t^1+{\frak g}\otimes t^4
+\cdots\\
&{\hat\frak g}_2= \cdots + {\frak g}\otimes t^{-1} +
{\frak g}\otimes t^2+{\frak g}\otimes t^5
+\cdots\endaligned$$

In terms of the affine root vectors discussed below (2.8),
 every element $E_{\alpha+m\delta}$ belongs
 to the subspace ${\hat\frak g}_{m\ mod\ N}$.

In the following section, I will illustrate the algorithms discussed
 here by using the specific algebra $\hat A_2$. I will also display
 a mixed grading of the {\it twisted} algebra $A_2^{(2)}$. In general,
 the underlying grading of a twisted algebra (which is neither vertical
 nor horizontal in the sense above)
 can be used to find
 graded contractions, but there are other gradings, which involve some
 mixing of this grading and other (vertical or horizontal) gradings.
 Apart from displaying an example at the end of the next section, I
 do not consider this problem any further here.

\medskip

\leftline {\it 2.3. An example: $\hat A_2$.}
\smallskip

A general element of the rank two, eight dimensional simple Lie algebra
 $A_2$ (or $\text {sl} (3,\Bbb C)$) can be written in the matrix form:
$$
\pmatrix
h_1 & e_{\alpha_1} & e_{\alpha_1 +\alpha_2}\\
e_{-\alpha_1} & h_2 & e_{\alpha_2}\\
e_{-(\alpha_1 +\alpha_2)} & e_{-\alpha_2} & -(h_1+h_2)
\endpmatrix.\tag2.13$$

Upon affinization, this algebra becomes the infinite dimensional
 Lie algebra
 $\hat A_2=\left (A_2\otimes\Bbb C\left [t,t^{-1}\right ]\right )
 \oplus \Bbb C k$,
 where the
 central element $k$ can be represented by the $3\times 3$ identity
 matrix. However the usual matrix product must be modified so as to
 satisfy the commutation relations (2.1).

The {\it horizontal} gradings of $\hat A_2$
 are provided by an EFO
 $[s_0, s_1, s_2]$, which describes a
 conjugacy class of elements of order $N=MC$, where
 $M=s_0+s_1+s_2$ and $C={3\over{gcd(3;s_1+2s_2)}}$ (see \cite{33}).
 The only element of order two is given by
 ${\bold s}=[0,1,1]$ and provides, according to (2.2), the grading:
$$
\aligned
&\left (A_2\right )_0\equiv {\frak h}+\Bbb C
 e_{\pm
 (\alpha_1+\alpha_2)},\\
&\left (A_2\right )_1\equiv \Bbb Ce_{\pm\alpha_1}+
 \Bbb Ce_{\pm\alpha_2}.\endaligned\tag2.14$$
$\frak h=\Bbb C h_1\dotplus\Bbb C h_2$ is the Cartan subalgebra of $A_2$.
In terms of the affine root vectors, the corresponding grading
 subspaces of $\hat\frak g$ are generated by:
$$
\aligned
&(\hat A_2)_0=\{E_{m\delta},
 E_{\pm(\alpha_1+\alpha_2 )+m\delta}, k\},\\
 &(\hat A_2)_1=\{E_{\pm\alpha_1+m\delta},
 E_{\pm\alpha_2+m\delta}\},\quad m\in\Bbb Z.\endaligned\tag2.14'$$
 Using $\delta=\alpha_0+\Psi_{A_2}=\alpha_0+\alpha_1+
 \alpha_2$, we can write, for instance,
 $E_{\alpha_1+m\delta}$ as
 $E_{m\alpha_0+(m+1)\alpha_1+m\alpha_2}$, etc.

The finite algebra $A_2$ admits two elements
 of order three, $[1,1,1]$ and
 $[0,1,0]$, which correspond to the decompositions
$$
\aligned
&\left (A_2\right )_0\equiv {\frak h},\\
&\left (A_2\right )_1\equiv \Bbb Ce_{\alpha_1}
 +\Bbb Ce_{\alpha_2}+\Bbb Ce_{-(\alpha_1+\alpha_2)},\\
&\left (A_2\right )_2\equiv \Bbb Ce_{-\alpha_1}
 +\Bbb Ce_{-\alpha_2}+\Bbb Ce_{\alpha_1+\alpha_2},
\endaligned\tag2.15$$
and
$$
\aligned
&\left (A_2\right )_0\equiv {\frak h}+
\Bbb C e_{\pm\alpha_2},\\
&\left (A_2\right )_1\equiv \Bbb Ce_{\alpha_1}
 +\Bbb Ce_{\alpha_1+\alpha_2},\\
&\left (A_2\right )_2\equiv \Bbb Ce_{-\alpha_1}
 +\Bbb Ce_{-(\alpha_1+\alpha_2)},\endaligned\tag2.16$$
respectively. Therefore, the $\Bbb Z_3$ grading of $\hat A_2$
 given by $[1,1,1]$ is
$$
\aligned
&(\hat A_2)_0=\{E_{m\delta},k\},\\
&(\hat A_2)_1=\{E_{\alpha_1+m\delta}, E_{\alpha_2+m\delta},
 E_{-\alpha_1-\alpha_2+m\delta}\},\\
&(\hat A_2)_2=\{E_{-\alpha_1+m\delta}, E_{-\alpha_2+m\delta},
 E_{\alpha_1+\alpha_2+m\delta}\},\quad m\in\Bbb Z,\endaligned\tag2.15'$$
and the grading
 $[0,1,0]$ is
$$
\aligned
&(\hat A_2)_0=\{E_{m\delta},
 E_{\pm\alpha_2+m\delta}, k\},\\
&(\hat A_2)_1=\{E_{\alpha_1
 +m\delta}, E_{\alpha_1+\alpha_2+m\delta}\},\\
&(\hat A_2)_2=\{E_{-\alpha_1+m\delta}, E_{-(\alpha_1+\alpha_2)
 +m\delta}\},\quad m\in\Bbb Z.\endaligned\tag2.16'$$
These expressions illustrate the fact that,
 for an horizontal grading, if $e_\alpha\in
{\frak g}_\mu$, then $E_{\alpha+m\delta}\in{\hat\frak g}_\mu$,
 for all $m\in\Bbb Z$.

{}From the section 2.2, the {\it vertical}
 $\Bbb Z_2$ grading is given by
$$
\aligned
&\left (\hat A_2\right )_0\equiv \left (A_2\otimes t^{2m}\right )
 \oplus \Bbb C k=\{E_{\alpha+2m\delta}, k\},\\
&\left (\hat A_2\right )_1\equiv A_2\otimes t^{2m+1}=
\{E_{\alpha+(2m+1)\delta}\},
\endaligned\tag2.17$$
the $\Bbb Z_3$ grading, by
$$
\aligned
&\left (\hat A_2\right )_0\equiv \left (A_2\otimes t^{3m}\right )
 \oplus \Bbb C k=\{E_{\alpha+3m\delta},k\},\\
&\left (\hat A_2\right )_1\equiv A_2\otimes t^{3m+1}=
 \{E_{\alpha+(3m+1)\delta}\},\\
&\left (\hat A_2\right )_2\equiv A_2\otimes t^{3m+2}=
 \{E_{\alpha+(3m+2)\delta}\},
\endaligned\tag2.18$$
and a general $\Bbb Z_N$ grading, by
$$
\aligned
&\left (\hat A_2\right )_0\equiv \left (A_2\otimes t^{Nm}\right )
 \oplus \Bbb C k=\{E_{\alpha+Nm\delta},k\},\\
&\left (\hat A_2\right )_1\equiv A_2\otimes t^{Nm+1}=
\{E_{\alpha+(Nm+1)\delta}\},\\
&\left (\hat A_2\right )_2\equiv A_2\otimes t^{Nm+2}=
\{E_{\alpha+(Nm+2)\delta}\},\\
& \ \vdots\quad\qquad\qquad\vdots\qquad\qquad\qquad\qquad\vdots\\
&\left (\hat A_2\right )_{N-1}\equiv A_2\otimes t^{Nm+N-1}=
\{E_{\alpha+(Nm+N-1)\delta}\}.\endaligned
\tag2.19$$

To illustrate the meaning of the expression (2.6), I now
 display a mixed
 vertical-horizontal grading. If we mix the horizontal decomposition
 (2.14') with the vertical grading (2.17),
 we get the following $\Bbb Z_2\otimes\Bbb Z_2$ grading:
$$
\aligned
&\left (\hat A_2\right )_{00}=\{E_{2m\delta}, E_{\pm (\alpha_1
 +\alpha_2)+2m\delta}, k\},\\
&\left (\hat A_2\right )_{01}=\{E_{(2m+1)\delta}, E_
{\pm (\alpha_1+\alpha_2)+(2m+1)\delta}\},\\
&\left (\hat A_2\right )_{10}=\{E_{\pm\alpha_1+2m\delta},
 E_{\pm\alpha_2+2m\delta}\},\\
&\left (\hat A_2\right )_{11}=\{E_{\pm\alpha_1+(2m+1)\delta},
 E_{\pm\alpha_2+(2m+1)\delta}\}; m\in\Bbb Z.\endaligned\tag2.20$$
The first $\Bbb Z_2$ index corresponds to the grading (2.14'), and
 the second index, to the decomposition (2.17).
 One can verify that (2.20) satisfies the relation (1.1).

I close this section by discussing some gradings of the twisted
 algebra $A_2^{(2)}$. The $\Bbb Z_2$ grading
 inherent to the twisting is
$$
\alignat{3}
& (A_2^{(2)})_0: && \quad t^j\otimes (h_{\alpha_1}+h_{\alpha_2}),
 t^j\otimes (e_{\alpha_1}+e_{\alpha_2}),
 (e_{-\alpha_1}+e_{-\alpha_2}),k, && j\ \text {even},\\
& (A_2^{(2)})_1: && \quad t^j\otimes (h_{\alpha_1}-h_{\alpha_2}),
 t^j\otimes (e_{\alpha_1}-e_{\alpha_2}),
 t^j\otimes (e_{-\alpha_1}-e_{-\alpha_2}), && \\
& && t^j\otimes e_{\alpha_1+\alpha_2},
 t^j\otimes e_{-(\alpha_1+\alpha_2)},&& j\ \text {odd}.\endalignat
$$

Oviously, this algebra admits other gradings. For instance, the grading
 above can be mixed with horizontal or vertical gradings.
 If we mix it
 with the $\Bbb Z_3$ grading (2.15), we get
 the $\Bbb Z_2\otimes\Bbb Z_3$ grading:
$$
\aligned
&(A_2^{(2)})_{00}=\Bbb C t^j\otimes (h_{\alpha_1}+h_{\alpha_2}),
 \quad j\ \text {even},\\
&(A_2^{(2)})_{01}=\Bbb C t^j\otimes (e_{\alpha_1}+e_{\alpha_2}),
 \quad j\ \text {even},\\
&(A_2^{(2)})_{02}=\Bbb C t^j\otimes (e_{-\alpha_1}+e_{-\alpha_2}),
 \quad j\ \text {even},\\
&(A_2^{(2)})_{10}=\Bbb C t^j\otimes (h_{\alpha_1}-h_{\alpha_2}),
 \quad j\ \text {odd},\\
&(A_2^{(2)})_{11}=\Bbb C t^j\otimes (e_{\alpha_1}-e_{\alpha_2})
+\Bbb C t^j\otimes e_{-(\alpha_1+\alpha_2)},
 \quad j\ \text {odd},\\
&(A_2^{(2)})_{12}=\Bbb C t^j\otimes (e_{-\alpha_1}-e_{-\alpha_2})
+\Bbb C t^j\otimes e_{\alpha_1+\alpha_2},
 \quad j\ \text {odd}.\endaligned
$$

\medskip

\leftline {\it 2.4. Gradings of representations.}
\smallskip

In this section I describe and give some examples of the
 gradings of integrable irreducible highest weight representations of
 untwisted affine Lie algebras, given a vertical or horizontal
 grading of the algebra.

An irreducible highest weight integrable module $V(\Lambda)$ of
 ${\hat\frak g}$ is labelled by its {\it highest weight}
 $\Lambda=(n;\Lambda_0,\dots,\Lambda_r )$, where
 $n,\Lambda_0,\dots,\Lambda_r$ are non-negative integers.
 As the finite case
 (see (2.3)), it can be weight decomposed as
$$
V(\Lambda)=\bigoplus_{\lambda\in{\hat\frak h}^*}V^\Lambda_\lambda,
\tag2.21$$
where the weight $\lambda=(n;\lambda_0,\dots,\lambda_r)$ has
 multiplicity $m^\Lambda_\lambda=$dim$V^\Lambda_\lambda$.
 ${\hat\frak h}$ is the Cartan subalgebra of ${\hat\frak g}$.
 An invariant
 of $V(\Lambda)$ is the {\it level}
 $\Lambda(k)=\sum_{j=0}^r\check c_j\Lambda_j$ ($k$: central
 element of ${\hat\frak g}$), where the $\check c_j$ are the {\it
 comarks} of ${\hat\frak g}$, defined by
 $\sum_{j=0}^rA_{jk}\check c_k=0$, for all $j$ ($A$: affine Cartan
 matrix of ${\hat\frak g}$). The integer $n$ in $\lambda$ is called
 {\it null depth}, and is equal to the number of $\alpha_0$'s that
 must be subtracted from $\Lambda$ to reach $\lambda$. The null depth
 determines the vertical gradings.

In order to grade the module $V(\Lambda)$ compatibly with some
 given grading of ${\hat\frak g}$, we consider the action of root
 vectors $E_\alpha$ on the vectors in $V(\Lambda)$ so as to coarsen
 the Cartan decomposition (2.21). To achieve this, we first express
 the roots in the basis of fundamental weights $\omega_0,\dots,
 \omega_r$:
$$
\aligned
&\alpha_0=(1;\alpha^0_0,\dots,\alpha^r_0)=
\delta+\alpha^0_0\omega_0+\cdots +\alpha^r_0\omega_r,\\
&\alpha_1=(0;\alpha^0_1,\dots,\alpha^r_1)=
\alpha^0_1\omega_0+\cdots +\alpha^r_1\omega_r,\\
&\quad\vdots\qquad\qquad\vdots\\
&\alpha_r=(0;\alpha^0_r,\dots,\alpha^r_r)=
\alpha^0_r\omega_0+\cdots +\alpha^r_r\omega_r,\endaligned
\tag2.22$$
where the coefficients are given by the affine Cartan matrix:
 $\alpha^j_k=A_{jk},$ for $j,k=0,\dots,r$, and
 $\delta =(1;0,\dots,0)$.
In the case of $\hat A_2$, the Cartan matrix is
$$
A=\pmatrix
 2 & -1 & -1\\
 -1 & 2 & -1\\
 -1 & -1 & 2\endpmatrix,\tag2.23$$
so that
$$
\aligned
&\alpha_0=(1;2,-1,-1),\\
&\alpha_1=(0;-1,2,-1),\\
&\alpha_2=(0;-1,-1,2).\endaligned\tag2.24$$

We use $E_\alpha$, expressed in the $\omega$-basis, and the fact
 that $E_\alpha\cdot V^\Lambda_\lambda\subseteq V^\Lambda_
{\lambda+\alpha}$ in order to find a compatible grading of $V(\Lambda)$.
For instance, the horizontal
 $\Bbb Z_2$ grading of $\hat A_2$ given by (2.14')
 reads, in this basis:
$$
\aligned
&(\hat A_2)_0=\{(m;0,0,0),(m;\mp 2,\pm 1,\pm 1),k\},\\
&(\hat A_2)_1=\{(m;\mp 1,\pm 2,\mp 1),(m;\mp 1,\mp 1,\pm 2)\},
\quad m\in\Bbb Z.\endaligned\tag2.25$$
In this simple case, we can find by inspection that
$$
\aligned
&V_0=\{(m;2\Bbb Z+1,\Bbb Z,\Bbb Z\},\\
&V_1=\{(m;2\Bbb Z,\Bbb Z,\Bbb Z\},\endaligned\tag2.26$$
together with (2.25) satisfies (1.2).

To illustrate the $\Bbb Z_3$ grading, we will consider the
 irrep $\Lambda=(1,0,0)$.
 Down to null depth $n=10$, its weight space decomposition
 has the form:

$$
\alignat{7}
& \lambda: &&\qquad w_0 &&\qquad w_1 &&\qquad w_2 &&
\qquad w_3 &&\qquad w_4 &&\qquad w_5\\
& n=0 &&\qquad 1 &&\qquad &&\qquad &&\qquad &&\qquad &&\qquad \\
& n=1 &&\qquad 2 &&\qquad 1 &&\qquad &&\qquad &&\qquad &&\qquad \\
& n=2 &&\qquad 5 &&\qquad 2 &&\qquad &&\qquad &&\qquad &&\qquad \\
& n=3 &&\qquad 10 &&\qquad 5 &&\qquad 1 &&\qquad &&\qquad &&\qquad \\
& n=4 &&\qquad 20 &&\qquad 10 &&\qquad 2 &&\qquad 1 &&\qquad &&\qquad \\
& n=5 &&\qquad 36 &&\qquad 20 &&\qquad 5 &&\qquad 2 &&\qquad &&\qquad \\
& n=6 &&\qquad 65 &&\qquad 36 &&\qquad 10 &&\qquad 5 &&\qquad &&\qquad \\
& n=7 &&\qquad 110 &&\qquad 65 &&\qquad 20 &&\qquad 10 &&\qquad 1 &&\qquad \\
& n=8 &&\qquad 185 &&\qquad 110 &&\qquad 36 &&\qquad 20 &&\qquad 2 &&\qquad \\
& n=9 &&\qquad 360 &&\qquad 185 &&\qquad 65 &&\qquad 36 &&\qquad 5 &&\qquad 1\\
& n=10 &&\qquad 481 &&\qquad 360 &&\qquad 110 &&\qquad 65 &&\qquad 10 &&\qquad
2
\endalignat
$$
where
$$
\aligned
w_0&= (1,0,0),\\
w_1&= (-1,1,1), (0,-1,2), (2,-2,1), (3,-1,-1),\\
w_2&= (-2,0,3), (1,-3,3), (4,0,-3),\\
w_3&= (-3,2,2), (-1,-2,4), (3,2,-4), (5,-2,-2),\\
w_4&= (-4,1,4), (-3,-1,5), (5,-5,1), (0,5,-4), (2,4,-5), (6,-1,-4),\\
w_5&= (-5,3,3), (-2,-3,6), (4,-6,3), (7,-3,-3),
\endaligned
$$
 and the vertical strings contain the weight multiplicities,
 given in \cite {29}.

The $\Bbb Z_3$ grading $[1,1,1]$ of (2.15') can be read
$$
\aligned
&(\hat A_2)_0=\{(m;0,0,0),k\},\\
&(\hat A_2)_1=\{(m;-1,2,-1),(m;-1,-1,2),(m;2,-1,-1)\},\\
&(\hat A_2)_2=\{(m;1,-2,1),(m;1,1,-2),(m;-2,1,1)\}.\endaligned
\tag2.27$$
To find the corresponding grading of $V(\Lambda)$,
 one can choose the highest weight
 $(0;1,0,0)$ to belong to $V_0$, and act iteratively on this weight
 with the elements of various grading subspaces (2.27). From (1.2),
 one finds, for all the weights down
 to $n=10$:
$$
\aligned
&V_0=\{(1,0,0),(-2,0,3),(1,3,-3),(4,0,-3),(-5,3,3),(-2,6,-3),\\
&\qquad (7,-3,-3),(4,3,-6),\dots\},\\
&V_1=\{(0,2,-1),(3,-1,-1),(-3,2,2),(3,2,-4),(-3,5,-1),(0,5,-4),\\
&\qquad(6,-1,-4),\dots\},\\
&V_2=\{(-1,1,1),(2,1,-2),(-1,4,-2),(-4,1,4),
(5,-2,-2),(2,4,-5),\\
&\qquad (5,1,-5),\dots\},\endaligned\tag2.28$$
plus all the permutations of the last two components $\lambda_1$ and
 $\lambda_2$ of each weight ({\it e.g.} $(4,-3,0)\in V_0$).
 The null depth $n$ is omitted in the weight because
 the grading does not depend on it. The grading has been chosen so
 that the weight $(1,0,0)$ belongs to $V_0$. The straightforward way to obtain
 (2.28) is by using (1.2) and apply all the elements of the different
 grading subspaces (${\hat\frak g}_1$ and
 ${\hat\frak g}_2$ from (2.27)) on $V_0$ so as to find a subset
 of each grading subspace of $V(\Lambda )$. Proceeding iteratively, we then
 apply the same elements (through (1.2)) on the identified elements
 of the $V_\mu$ found in the first step, to find further elements of $V_\mu$.
 The grading (2.28) lies in the direction of the
 {\it principal slicing} \cite{29} of the weight
 system.

For the grading $[0,1,0]$ of (2.16') we have
$$
\aligned
&(\hat A_2)_0=\{(m;0,0,0),(m;\pm 1,\pm 1,\mp 2),k\},\\
&(\hat A_2)_1=\{(m;-1,2,-1),(m;-2,1,1)\},\\
&(\hat A_2)_2=\{(m;1,-2,1),(m;2,-1,-1)\}.\endaligned
\tag2.29$$
By proceeding as for (2.28), we find the decomposition:
$$
\aligned
&V_0=\{(1,0,0),(0,-1,2),(2,1,-2),(-1,-2,4),(3,2,-4),(-4,4,1),\\
&\qquad (-3,5,-1),(5,-5,1),(6,-4,-1),(-5,3,3),(-2,-3,6),(-2,6,-3),\\
&\qquad (4,3,-6),(4,-6,3),(7,-3,-3),\dots\},\\
&V_1=\{(-1,1,1),(0,2,-1),(-2,0,3),(1,3,-3),(4,-3,0),(3,-4,2),\\
&\qquad (5,-2,-2),(-3,-1,5),(2,4,-5),(2,-5,4),(6,-1,-4),\dots\},\\
&V_2=\{(2,-2,1),(3,-1,-1),(-2,3,0),(-3,2,2),(1,-3,3),(-1,4,-2),\\
&\qquad (4,0,-3),(-4,1,4),(0,-4,5),(0,5,-4),(5,1,-5),\dots\}.\endaligned
\tag2.30$$
Again, the grading is independent of the null depth.

{}From (2.19), we see that a general vertical
 $\Bbb Z_N$ grading has the form
$$
\aligned
&(\hat A_2)_0=\{(Nm,\alpha)\},\\
&(\hat A_2)_1=\{(Nm+1,\alpha)\},\\
&\quad\vdots\qquad\vdots\\
&(\hat A_2)_k=\{(Nm+k,\alpha)\},\endaligned\tag2.31$$
where $m\in\Bbb Z$, and $\alpha$ is any root of $\hat A_2$.
 Now the grading depends on the null depth only.
 The corresponding grading of $V(\Lambda )$ is
$$
V_k=\{(k,\alpha)\},\qquad k=0,\dots,N-1\ mod\ N, \text {for all}\ \alpha.
\tag2.32$$

\bigskip

\leftline {\bf 3. GRADED CONTRACTIONS.}
\medskip

In section 1, I have defined the graded contractions of any
 Lie algebra and its representations. In section 2, I have described
 the horizontal and the vertical gradings, and,
 more particularly, the
 toroidal gradings of Kac-Moody algebras and their
 irreducible representations. These are the basic ingredient needed
 to contract an algebra and its representations.
 It is now straightforward to obtain the
 graded contractions of Kac-Moody algebras, which
 form a new class of infinite dimensional Lie algebras.
 To summarize the contraction of algebras: one gets an horizontal grading
 (2.8) by using the expression (2.2) to find the eigenspaces of the EFO,
 or a vertical grading by using (2.11)
 and (2.8'). To find the graded contractions, one just replace the
 solutions of (1.5) in the modified commutation relations (1.3). The grading of
 representations has been described and illustrated in section 2.4.

The purpose of this section is {\it not} to display huge lists of
 contractions, but rather to describe their general properties. There
 exists a computer program \cite{24} that provides the solutions of
 equations (1.5), (1.6) and (1.9), given the grading group $\Gamma$ and
 the structure of the grading ({\it i.e.} generic or non-generic).
 Each solution then provides a contraction of the algebra
 or the representation.

The most straightforward definition of graded contractions of Kac-Moody
 algebras is, after (1.3), to deform the commutator (2.1) as
$$
\aligned
 [a\otimes t^m, b\otimes t^n]_\varepsilon
&\equiv \varepsilon_{\mu ,\nu}[a\otimes t^m, b\otimes t^n]\\
&=\varepsilon_{\mu ,\nu}
 [a,b]\otimes t^{m+n}+\varepsilon_{\mu ,\nu}mkB(a,b)\delta_{m+n,0},
\endaligned
 \tag3.1$$
 where $a\otimes t^m\in{\frak g}_\mu, b\otimes t^n\in{\frak g}_\nu$,
 a vertical or horizontal grading.
 (As discussed below, an interesting alternative is to deform
 simultaneously the Killing form $B$\cite {32},
 so that it becomes possible to preserve the central extension
 whereas the first term in (3.1) is put to zero).

\medskip
\leftline {\it Comparison with In\"on\"u-Wigner contractions.}
\smallskip

First, we compare the In\"on\"u-Wigner contraction of a Kac-Moody
 algebra (studied in \cite{16}) with the particular case
 of a $\Bbb Z_2$ graded contraction. We write the basis of the algebra
 ${\hat\frak g}$ as $T^a_m$,
 where $\ a=1,\dots , \text {dim} {\hat\frak g},\
 \text {and}\ m\in\Bbb Z$, and with the commutation relations:
$$
[T^a_m,T^b_n]=if^{a,b}_cT^c_{m+n}+\frac 12km\delta^{a,b}
\delta_{m+n}.\tag3.2$$
Then, we decompose ${\hat\frak g}$ {\it \`a la} In\"on\"u-Wigner,
 by writing the underlying vector space as ${\hat\frak g}=
 {\hat\frak g}_0+{\hat\frak g}_1$, where
 ${\hat\frak g}_0=\{T^\alpha_m\}, \alpha =1,2,\dots ,r$, forms a
 subalgebra of ${\hat\frak g}$, and
 ${\hat\frak g}_1=\{T^i_m\}, i=r+1, r+2,\dots , \text {dim} {\hat\frak g}$
 is its complementary subspace.
 The commutation relations (3.2) must take the form
$$
\aligned
&[T^\alpha_m, T^\beta_n]=if^{\alpha,\beta}_\gamma T^\gamma_{m+n}
 +\frac 12km\delta^{\alpha,\beta}\delta_{m+n,0},\\
&[T^\alpha_m,T^i_n]=if^{\alpha,i}_jT^j_{m+n},\\
&[T^i_m, T^j_n]=if^{i,j}_\alpha T^\alpha_{m+n}+\frac 12km
\delta^{i,j}\delta_{m+n,0},\endaligned\tag3.3$$
in order to define an In\"on\"u-Wigner contraction in this basis.
 The contraction is then defined by multiplying all the basis elements
 of the vector subspace ${\hat\frak g}_1$
 by a contraction parameter
 $\varepsilon$, and, in the limit $\varepsilon\rightarrow 0$, the
 commutators in the third row of (3.3) vanish. Thus, the resulting
 algebra admits a $\Bbb Z_2$ grading, where
 ${\hat\frak g}_0=\{T^\alpha_m, k\}$
 and ${\hat\frak g}_1=\{T^i_n\}$.
 Therefore, one can say that the In\"on\"u-Wigner
 contraction of an affine algebra is a particular case of $\Bbb Z_2$
 graded contraction, with
 $\varepsilon_{0,0}=1=\varepsilon_{0,1}$ and
 $\varepsilon_{1,1}=0$.  Obviously, there are other graded contractions
 which lead to an In\"on\"u-Wigner contraction.

As mentioned in \cite {25} for the particular case of
 $\hat A_1$, among the contractions of a general algebra ${\hat\frak g}$,
 there are semi-direct products of the initial ${\hat\frak g}$
 with an infinite dimensional abelian ideal, or ``translation''
 algebra. In order words, among the possible contractions of
 ${\hat\frak g}$, one finds the algebra ${\hat\frak g}\triangleright
 {\frak a}$, where ${\frak a}$ is an infinite dimensional abelian
 ideal of the contracted algebra. This may be surprising because
 it is specific to the infinite dimensional Lie algebras, and cannot
 occur in the finite cases.
 This occurs when we take a {\it vertical} grading, where, from (2.8),
 ${\hat\frak g}_0=({\frak g}\otimes\Bbb C[t,t^{-1}]_0)\oplus\Bbb Ck
 \approx {\hat\frak g}$. If
 $\varepsilon_{0,\mu}=1$ for all $\mu$, and all the other
 parameters vanish, then
 the subalgebra ${\hat\frak g}$ will be preserved, and so do the
 commutators involving this subalgebra and the remaining basis elements.
 Since the remaining commutators all vanish, the corresponding ideal
 is abelian.

In fact, the graded contractions allow to go much further than this.
 Whenever $\varepsilon_{0,0}=1$, the subalgebra ${\hat\frak g}$
 ({\it i.e.} the original algebra) will be contained in the contracted
 algebra, either in direct or semi-direct sums.
 I will illustrate this
 with $\Bbb Z_2$ contractions. In addition to the contraction
 mentioned previously, there are two other non-trivial contractions,
 namely one where $\varepsilon_{0,0}=1, \varepsilon_{0,1}=0=\varepsilon_{1,1}$,
 and the other, with $\varepsilon_{0,0}=0=\varepsilon_{0,1},
 \varepsilon_{1,1}=1$. For the first, only the subalgebra ${\hat\frak g}$
 is preserved, so that the contracted algebra is
 ${\hat\frak g}^\varepsilon ={\hat\frak g}\oplus{\frak a}$, where
 ${\frak a}$ is abelian (and, obviously, infinite). Under this contraction
 the vector space underlying ${\frak a}$ carries no longer a representation
 space of the subalgebra ${\hat\frak g}$ in the adjoint representation.
 In the second case, the only
 commutation relations that are not deformed to zero are $[{\hat\frak g}_1,
 {\hat\frak g}_1]$.
 Therefore, the subspace
 ${\hat\frak g}_0$ becomes abelian, and the commutation relations involving
 any of its elements also vanish.

We note also that the center is modified under a contraction,
 as in the finite case.
 Whereas the center of an affine algebra consists only of its
 central extension, it usually becomes
 bigger after a contraction. For instance,
 in the first $\Bbb Z_2$ contraction displayed in the previous paragraph
 the center also includes all the elements of the subspace ${\hat\frak g}_1$.
 In the second $\Bbb Z_2$ contraction, the center includes the
 subalgebra ${\hat\frak g}_0$ ({\it i.e.} the original algebra).
 Depending on the particular grading that is preserved, and depending on
 the original algebra, there might be additional elements in
 the center.

\medskip
\leftline {\it Generators of positive root vectors.}
\smallskip

Another interesting property of a Kac-Moody algebra that is modified
 under a contraction is the minimal set of generators of positive root
 vectors. In general, the greater the number of contracted
 commutators ({\it i.e.} zero after contraction),
 the greater is the set of
 such generators. Below, I illustrate this point by discussing
 in detail some examples with
 ${\hat A}_1$ and ${\hat A}_2$.

The set of positive root vectors of ${\hat A}_1$ is given by
 $E_{p\alpha_0+q\alpha_1}$, where $p$ and $q$ are positive integers
 such that $-1\leq p-q\leq 1$. The gradings of ${\hat A}_1$ in terms
 of these vectors are easy to visualize if we write them as:
$$
\alignat {6}
& E_{\alpha_1} && \quad E_{\alpha_1+\delta} && \quad E_{\alpha_1+2\delta}
 && \quad E_{\alpha_1+3\delta} && \quad E_{\alpha_1+4\delta} &&
 \quad \cdots \\
&  && \quad E_{\delta} && \quad E_{2\delta}
 && \quad E_{3\delta} && \quad E_{4\delta} &&
 \quad \cdots \\
&  && \quad E_{-\alpha_1+\delta} && \quad E_{-\alpha_1+2\delta}
 && \quad E_{-\alpha_1+3\delta} && \quad E_{-\alpha_1+4\delta} &&
 \quad \cdots
\endalignat
$$

Now I will show explicitly how to obtain the generators for all the
 $\Bbb Z_2$ contractions. (This was done in \cite {25} but
 I obtain it here in a more sytematical way, easy to generalize to
 other algebras and gradings.)

For the {\it horizontal} grading, ${\hat\frak g}_0=\{E_{m\delta},
 k;\ m\geq 1\}$
 are the elements of the second row,
 and ${\hat\frak g}_1=\{E_\alpha , E_{\alpha +m\delta},
 E_{-\alpha+m\delta};\ m\geq 1\}$ corresponds to
 the first and third rows.
 To find the generators of positive root vectors
 for the contraction $\varepsilon_{0,0}=1=\varepsilon_{0,1},\
 \varepsilon_{1,1}=0$, we consider each element of the array, one at the
 time, and see if it can be obtained by commutation of previous
 generators, by taking into account the contraction parameters. It is
 convenient to start from the left, and from the bottom to the top.
 So the two elements that we first keep are $E_{\alpha_1}$ and
 $E_{-\alpha_1+\delta}$. The next two elements, $E_\delta$ and
 $E_{\alpha_1+\delta}$ can be obtained by the commutator of the
 first two, so we do not keed them as generators. The
 next (and the last) one to be retained is $E_{-\alpha_1+\delta}$,
 which cannot be obtained from any commutator of the other elements.
 Therefore, the set of positive roots of the graded contractions
 is generated by three vectors: $E_{\alpha_1},\
 E_{-\alpha_1+\delta}$, and $E_{-\alpha_1+2\delta}$.

By a similar reasoning, we find that the generators corresponding
 to the contraction $\varepsilon_{0,0}=1, \varepsilon_{0,1}=0=
\varepsilon_{1,1}$ are
$$
\aligned
& E_{\alpha_1},\  E_{\alpha_1+(2k+1)\delta},\ E_{-\alpha_1+2\delta},\\
  & E_{(2k+1)\delta},\ E_{-\alpha_1+(2k+1)\delta},
\quad k\geq 0.\endaligned
$$
For the contraction  $\varepsilon_{0,0}=0=\varepsilon_{0,1}$ and
 $\varepsilon_{1,1}=1$, the generators are
$$
E_{\alpha_1},\ E_{\alpha_1+(2k+1)\delta},
\ E_{(2k+1)\delta},\ E_{-\alpha_1+(2k+1)\delta},\
\quad k\geq 0.
$$

The {\it vertical} grading, for which ${\hat\frak g}_0=\{E_{\alpha_1},
 E_{\alpha_1+2m\delta}, E_{2m\delta}, E_{-\alpha_1+2m\delta},
 k;\ m\geq 1\}$ and ${\hat\frak g}_1=\{E_{\alpha_1+(2m-1)\delta},
 E_{(2m-1)\delta}, E_{-\alpha_1+(2m-1)\delta};\ m\geq 1\}$
 consists in the elements of the odd columns
 for ${\hat\frak g}_0$, and the even columns for
 ${\hat\frak g}_1$.
 For this grading (which is non-generic, because $[{\hat\frak g}_0,
 {\hat\frak g}_0]=0$), we have two contractions
 $\pmatrix
 \varnothing && \varepsilon_{0,1}\\
 \varepsilon_{0,1} && \varepsilon_{1,1}
\endpmatrix$,
 together with
 their respective generators of positive root vectors:
$$
\alignat{2}
&\pmatrix
 \varnothing && 1\\
 1 && 0\endpmatrix\quad \text {with}
\quad &&
E_{\alpha_1},\ E_{(k+1)\delta},\ E_{-\alpha_1+\delta},\\
&\pmatrix
 \varnothing && 0\\
 0 && 1\endpmatrix\quad \text {with}
\quad &&
E_{\alpha_1+k\delta},\ E_{-\alpha_1+(k+1)\delta},
\qquad k\geq 0.
\endalignat
$$

I now illustrate this with the algebra
 ${\hat A}_2$, and its $\Bbb Z_2$ gradings (2.14) and (2.17).
 As before, it is convenient to display the positive root
 vectors of ${\hat A}_2$ in an array:
$$
\alignat{6}
& E_{\alpha_1+\alpha_2} && E_{\alpha_1+\alpha_2+\delta} &&
 E_{\alpha_1+\alpha_2+2\delta} && E_{\alpha_1+\alpha_2+3\delta}
 && E_{\alpha_1+\alpha_2+4\delta} && \cdots\\
& E_{\alpha_2} && E_{\alpha_2+\delta} && E_{\alpha_2+2\delta} &&
 E_{\alpha_2+3\delta} && E_{\alpha_2+4\delta} && \cdots\\
& E_{\alpha_1} && E_{\alpha_1+\delta} && E_{\alpha_1+2\delta} &&
 E_{\alpha_1+3\delta} && E_{\alpha_1+4\delta} && \cdots\\
&  && E_{\delta} && E_{2\delta} &&
 E_{3\delta} && E_{4\delta} && \cdots\\
&  && E_{-\alpha_1+\delta} && E_{-\alpha_1+2\delta} &&
 E_{-\alpha_1+3\delta} && E_{-\alpha_1+4\delta} && \cdots\\
&  && E_{-\alpha_2+\delta} && E_{-\alpha_2+2\delta} &&
 E_{-\alpha_2+3\delta} && E_{-\alpha_2+4\delta} && \cdots\\
&  && E_{-(\alpha_1+\alpha_2)+\delta} &&
 E_{-(\alpha_1+\alpha_2)+2\delta} && E_{-(\alpha_1+\alpha_2)+3\delta}
 && E_{-(\alpha_1+\alpha_2)+4\delta} && \cdots\endalignat
$$

The {\it horizontal}
 grading provided by the EFO $[0,1,1]$ (see (2.14'))
 consists in the top, middle and bottom rows for ${\hat\frak g}_0$,
 and the four remaining rows for ${\hat\frak g}_1$. It is a generic
 grading for which the generators are:
$$
E_{\alpha_1},\ E_{\alpha_2},\ E_{\alpha_1+\alpha_2},
 E_{-(\alpha_1+\alpha_2)+\delta},
$$
for the contraction
$\pmatrix
 \varepsilon_{0,0} && \varepsilon_{0,1}\\
 \varepsilon_{0,1} && \varepsilon_{1,1}\endpmatrix
=
\pmatrix
 1 && 1\\
 1 && 0\endpmatrix$,
$$
\aligned
& E_{\alpha_1},\ E_{\alpha_2},\ E_{\alpha_1+\alpha_2},
 E_{-(\alpha_1+\alpha_2)+\delta},\\
& E_{\alpha_1+k\delta},\ E_{\alpha_2+k\delta},\
 E_{-\alpha_1+k\delta},\ E_{-\alpha_2+k\delta},
\quad k\geq 1,\endaligned
$$
for the contraction
$\pmatrix
 1 && 0\\
 0 && 0\endpmatrix$, and
$$
\aligned
& E_{\alpha_1},\ E_{\alpha_2},
\ E_{-(\alpha_1+\alpha_2)+\delta},\\
& E_{\alpha_1+k\delta},\ E_{\alpha_2+k\delta},\
 E_{-\alpha_1+k\delta},\ E_{-\alpha_2+k\delta},
 \qquad k\geq 1,\endaligned
$$
for the contraction
$\pmatrix
 0 && 0\\
 0 && 1\endpmatrix$.

The {\it vertical} $\Bbb Z_2$ grading (2.17) has the elements
 of ${\hat\frak g}_0$ given by the odd columns, and the elements
 of ${\hat\frak g}_1$ given by the even columns.
 It is another generic
 grading for which the generators are:
$$
E_{\alpha_1},\ E_{\alpha_2},\ E_{-(\alpha_1+\alpha_2)+\delta},
 \ E_\delta,\ E_{2\delta},\ E_{-(\alpha_1+\alpha_2)+2\delta},
$$
for the contraction
$\pmatrix
 1 && 1\\
 1 && 0\endpmatrix$,
$$
E_{\alpha_1},\ E_{\alpha_2},\ E_{-(\alpha_1+\alpha_2)+2\delta},
 \ E_{\Delta+(2k-1)\delta},\qquad k\geq 1,
$$
where $\Delta$ represents all the roots (including the zero root)
 of $\frak g$, for the contraction
$\pmatrix
 1 && 0\\
 0 && 0\endpmatrix$, and
$
E_{\Delta+2k\delta},\qquad k\geq 1,
$
for the contraction
$\pmatrix
 0 && 0\\
 0 && 1\endpmatrix$.

\medskip
\leftline {\it Contractions of extended algebras.}
\smallskip

An interesting aspect of the contraction of affine Kac-Moody algebras
 is the behaviour, under a contraction, of the {\it extended} Kac-Moody algebra
 ${\hat\frak g}^e={\hat\frak g}\triangleleft\frak V$, where
 $\frak V$ is the Virasoro algebra associated to ${\hat\frak g}$.
 From the structure of semi-direct sum, we can
 see that, given a {\it vertical} $\Bbb Z_N$ grading (2.8), the corresponding
 $\Bbb Z_N$ grading of ${\hat\frak g}^e$ is given by
$$
\aligned
&{\hat\frak g}^e_0=\left ({\frak g}_0\otimes
 \Bbb C\left [t,t^{-1}\right ]\right )\oplus \Bbb Ck
\oplus \Bbb CL_{0\ mod\ N},\\
&{\hat\frak g}^e_{j\neq 0}={\frak g}_j\otimes
 \Bbb C\left [t,t^{-1}\right ]\oplus \Bbb C
 L_{j\ mod\ N},\qquad j=1,\dots,N-1,\endaligned\tag2.8$$
On the other hand, in the case of an horizontal grading, the
 Virasoro algebra $\frak V$ is contained completely in the
 grading subspace ${\hat\frak g}_0^e$.

It would be very interesting to study the representations of
 $\frak V$ obtained through the Sugawara construction, and see if the
 conclusions above are manifest from this construction. Let us
 just show how such a study should proceed.

The basis elements $L_n$ of $\frak V$ are defined by
$$
L_n=\frac 12\sum_{m\in\Bbb Z} :
 \sum_{i, j=1}^{\text {dim}\frak g} a_i\otimes t^{-m}
 a_j\otimes t^{m+n}:B(a_i, a_j),$$
where $B(\cdot ,\cdot)$ is the Killing form of ${\frak g}$,
 which is just the kronecker symbol in the bases usually utilized
 in the Sugawara construction. However, because a general grading
 of ${\frak g}$ is not always associated to such a basis, we must
 keep it explicitly in the sum.

Next we have to examine the steps of this construction, by taking
 into account the contraction parameters introduced both in the
 commutation relations (1.3), the action on the representation (1.4), and
 the bilinear form (see \cite {32}).
 Each term in the sum then takes the form:
$$\sigma^\psi(:{\frak g}_\mu{\frak g}_\nu:)V_\rho=\varepsilon_{\mu,\nu}
\psi_{\mu+\nu,\rho}\gamma_{\mu,\nu}
\sigma(:{\frak g}_\mu{\frak g}_\nu:)V_\rho.$$
The detailed investigation of this construction is beyond the scope
 of the present paper. I plan to study the contraction of Sugawara
 (studied with the traditional method in \cite {19}), and related
 constructions (GKO, Virasoro) soon.

In relation with the construction of WZW models,
 there are strong indications that the family of solvable Lie algebras
 introduced in \cite {34} can be obtained through a graded
 contraction, although they cannot be obtained from a In\"on\"u-Wigner
 contraction.
 For instance, a $\Bbb Z_3$ graded structure is
 inherent in these algebras. However, for the algebras $\Cal A_{3m},\
 (m\geq 3)$ the possible algebras to be contracted (which have the
 correct dimension) is huge, and there is no systematical way to
 identify them.
 This study
 is postponed to a future work.

\medskip
\leftline {\it Contractions involving a deformation of the bilinear form.}
\smallskip

As mentioned at the beginning of this section, we can define the
 contracted commutators by allowing the invariant bilinear form
 $B$ to be contracted as well \cite{32}. Given
 an horizontal grading, with $a\in{\frak g}_\mu$ and
 $b\in{\frak g}_\nu$,
 the commutator (3.1) is then modified to
$$
[a\otimes t^m, b\otimes t^n]_\varepsilon=\varepsilon_{\mu ,\nu}
 [a,b]\otimes t^{m+n}+\varepsilon_{\mu ,\nu}
\gamma_{\mu ,\nu}mkB(a,b)\delta_{m+n,0},
\tag3.4$$
where $B$ is replaced by $B^\gamma=\gamma B$. This permits to preserve
 the second term on the right-hand side even if $\varepsilon =0$,
 by choosing $\gamma$ such that $\gamma\varepsilon$ is constant.
 From \cite{32}, $B^\gamma({\frak g}_\mu,{\frak g}_\nu )
\equiv\gamma_{\mu ,\nu}B({\frak g}_\mu,{\frak g}_\nu )$, where
 $\gamma$ must satisfy \cite{32}
$$
\aligned
&\varepsilon_{\mu,\nu}\gamma_{\mu+\nu,\sigma}=
\varepsilon_{\nu,\sigma}\gamma_{\mu,\nu+\sigma},\\
&\gamma_{\mu,\nu}=\gamma_{\nu,\mu}.\endaligned
\tag3.5$$

To illustrate this, consider the $\Bbb Z_2$ contraction:
 $\varepsilon_{0,0}=1, \varepsilon_{0,1}=0=\varepsilon_{1,1}.$
 The corresponding solutions of (3.5) for $\gamma$ are: $\gamma_{0,1}=0$,
 with $\gamma_{0,0}, \gamma_{1,1}$ free. One can choose
 $\gamma_{1,1}$ to approach the infinity as $\varepsilon_{1,1}$
 approaches zero, such that $\varepsilon_{1,1}\gamma_{1,1}=K$,
 a constant.
 The commutators of the contracted algebra then become
$$
\aligned
&[(a\otimes t^m)_0,(b\otimes t^n)_0]_\varepsilon =
 [a,b]\otimes t^{m+n}+\gamma_{0,0}mkB(a,b)\delta_{m+n,0},\\
&[(a\otimes t^m)_0,(b\otimes t^n)_1]_\varepsilon =0,\\
&[(a\otimes t^m)_1,(b\otimes t^n)_1]_\varepsilon =KmkB(a,b)
 \delta_{m+n,0}.\endaligned\tag3.6$$

The oscillator algebra, as mentioned at the very
 beginning of this paper, can be obtained through the trivial
 $\Bbb Z_2$ contraction $\varepsilon_{0,0}=\varepsilon_{0,1}=
\varepsilon_{1,1}=0$ (for which the three parameters $\gamma$
 are free) with $\gamma_{0,1}=0=\gamma_{1,1}$ and $\gamma_{0,0}$
 approaches the infinity as $\varepsilon_{0,0}$ approaches $0$, so
 that $\varepsilon_{0,0}\gamma_{0,0}=K$. The central term above
 is then preserved. Actually, when done that way the oscillator
 algebra is a subalgebra of the initial Kac-Moody algebra. To get
 the true oscillator algebra, one just takes the trivial grading
 ${\hat\frak g}_0={\hat\frak g},\ {\hat\frak g}_1=0$, which shows
 how even this important algebra is a rather trivial
 graded contraction.

\enddocument

\Refs

\ref\no1\by E. In\"on\"u and E. P. Wigner\paper On the contraction
 of groups and their representations
\jour Proc. Nat. Acad. Sci. US\vol 39\yr 1953
\pages 510--524\endref

\ref\no2\by R. Gilmore\paper Lie Groups, Lie Algebras and Some of Their
 Applications\jour (Chap. 10), J. Wiley \& Sons, New York,\yr 1974\endref

\ref\no3\by E. J. Saletan\paper Contractions of Lie groups\jour
 J. Math. Phys.\yr 1961\vol 2\pages 1-21
 \endref

\ref\no4\by H. D. Doebner and O. Melsheimer\paper On a class of
 generalized group contractions\jour Nuov. Cim.\vol A 49\yr 1967
 \pages 306--311\endref

\ref\no5\by G. C. Hegerfeldt\paper Some properties of a class of
 generalized In\"on\"u-Wigner contractions\jour Nuov. Cim.
\vol A 51\yr 1967\pages 439--447\endref

\ref\no6\by M. L\'evy-Nahas\paper Deformation and contraction of
 Lie algebras\jour J. Math. Phys.\vol 8\yr 1967\pages 1211--1222
 \endref

\ref\no7\by E. Celeghini and M. Tarlini\paper Contraction of group
 representations I\jour Nuov. Cim.\vol B 61\yr 1981\pages 265--277
 \moreref\paper II,\vol B 65\yr 1981\pages 172--180\moreref\paper III,
 \vol B 68\yr 1982\pages 133--141\endref

\ref\no8\by H. Bacry and J. M. L\'evy-Leblond\paper Possible
 kinematics\jour J. Math. Phys.\vol 9\yr 1968
 \pages 1605--1614\endref

\ref\no9\by S. Str\"om\paper Construction of representations of the
 inhomogeneous Lorentz group by means of contraction of representations
 of the $(1+4)$ de Sitter group
 \jour Ark. Fys.
 \vol 30\yr 1965\pages 455--472\endref

\ref\no10\by R. Mignani\paper Instability of invariance groups of
 space-time, group contractions and models of universe\jour
 Nuov. Cim. Lett.\vol 23\yr 1978\pages 349--352\endref

\ref\no11\by A. Bohm and R.R. Aldinger\paper Examples of group
 contraction\jour Proc. of 11th Int. Coll. Group Theor. Meth.
 Phys., Instanbul\yr 1982\endref

\ref\no12\by E. Celeghini, P. Magnollay, M. Tarlini and
 G. Vitiello\paper Nonlinear realization and contraction of
 group representations\jour Phys. Lett.\vol B 162\yr 1985
 \pages 133--136\endref

\ref\no13\by L. D. Faddeev and N. Yu.
 Reshetikhin\paper Integrability of the principal chiral
 field model in (1+1)-dimension\jour Ann. Phys.
 \vol 167\yr 1986\pages 227--256\endref

\ref\no14\by N. A. Gromov\paper Contractions and analytic
 continuations of classical groups (in Russian), Syktyvkar
 \yr 1990\endref

\ref\no15\by M. de Montigny\paper The de Sitter-invariant
 differential equations and their contraction to Poincar\'e
 and Galilei\jour Nuov. Cim.\vol B 108\yr 1993
 \pages 1171--1180\endref

\ref\no16\by P. Majumdar\paper In\"on\"u-Wigner contraction of
 Kac-Moody algebras
 \jour J. Math. Phys.\vol 34\yr 1993
 \pages 2059--2065\endref

\ref\no17\by E. Witten\paper Non-abelian bosonization in two dimensions
 \jour Comm. Math. Phys.\vol 92\yr 1984\pages
 455--472\endref

\ref\no18\by C. R. Nappi and E. Witten\paper A WZW model based on a non
 semi-simple group\jour Phys. Rev. Lett.\vol 71\yr 1993
 \pages 3751--3753\endref

\ref\no19\by D. I. Olive, E. Rabinovici and A. Schwimmer
 \paper A class of string backgrounds as a semiclassical limit
 of WZW models\jour Phys. Lett. B\vol 321\yr 1994\pages 361--364
 \endref

\ref\no20\by  K. Sfetsos\paper Exact string backgrounds from WZW
 models based on non-semi-simple group\jour Int. J. Mod. Phys. A
 \vol 9\yr 1994\pages 4759--4766
 \moreref \paper Gauged WZW models and
 non-abelian duality\jour Phys. Rev. D\vol 50\yr 1994\pages
2784--2798
\moreref \paper Gauging a non-semisimple WZW model\jour Phys. Lett. B
\vol 324\yr 1994\pages 335--344\endref

\ref\no21\by J. Patera and H. Zassenhaus\paper On Lie gradings I
 \jour Lin. Alg. Appl.\vol 112\yr 1989\pages 87--159\endref

\ref\no22\by M. de Montigny and J. Patera\paper Discrete and continuous
 graded contractions of Lie algebras and superalgebras
\jour J. Phys. A: Math. Gen.
 \vol 24\yr 1991\pages 525--547\endref

\ref\no23\by R. V. Moody and J. Patera\paper Discrete and continous
 graded contractions of representations of Lie algebras
\jour J. Phys. A: Math. Gen.
 \vol 24\yr 1991\pages 2227--2257\endref

\ref\no24\by D. B\'erub\'e and M. de Montigny\paper The computer
 calculation of graded contractions of Lie algebras and their
 representations\jour Comp. Phys. Comm.\vol 76\yr 1993
 \pages 389--410\endref

\ref\no25\by A. Hussin, R. C. King, X. Leng and J. Patera\paper
 Graded contractions of the affine Lie algebra $A_1^{(1)}$, it
 representations and tensor products, and an application to the
 branching rule $A_1^{(1)}\subseteq A_1^{(1)}$\jour J. Phys. A:
 Math. Gen.\vol 27\yr 1994\pages 4125--4152\endref

\ref\no26\by V. G. Kac\paper Infinite dimensional Lie algebras\jour
 Cambridge Univ. Press, Cambridge\yr 1990\endref

\ref\no27\by P. Goddard and D. Olive\paper Kac-Moody and Virasoro
 algebras in relation to quantum physics\jour Int. J. Mod.
 Phys.\vol A1\yr 1986\pages 303--414\endref

\ref\no28\by J. F. Cornwell\paper Group Theory in Physics, Vol. 3
 \jour New York\yr 1989\endref

\ref\no29\by S. Kass, R. V. Moody, J. Patera and R. Slansky\paper
 Affine Lie algebras, weight multiplicities and branching rules,
 Vol. I \& II\jour Univ. of California Press, Berkeley\yr 1990\endref

\ref\no30\by V. Kac\paper Automorphisms of finite order of semisimple
 Lie algebras\jour J. Func. Anal. Appl.\yr 1969\pages 252-254\endref

\ref\no31\by J. Patera and R. V. Moody\paper Characters of elements
 of finite order in Lie groups\jour SIAM J. Alg. Disc. Meth.\vol 5\yr
 1984\pages 359--383\endref

\ref\no32\by M. de Montigny\paper Graded contractions of bilinear
 invariant forms of Lie algebras\jour J. Phys. A: Math. Gen.
 \vol 27\yr 1994\pages 4537--4548\endref

\ref\no33\by M. de Montigny\paper Elements of finite order in Lie groups
 and discrete gauge symmetries\jour Nucl. Phys.
 \vol B 439\yr 1995\pages 665--676\endref

\ref\no34\by A. Giveon, O. Pelc and E. Rabinovici\paper WZNW models
 and gauged WZNW models based on a family of solvable Lie algebras
\jour Racah Institute, RI-8-95, hep-th/ 9509013\yr 1995\endref

\endRefs
\bye